\documentclass[12pt]{article}
\usepackage{amssymb,amsmath}

\newcommand{\Q}{\mathbb{Q}} %
\newcommand{\R}{\mathbb{R}} %
\newcommand{\C}{\mathbb{C}} %
\newcommand{\N}{\mathbb{N}} %
\renewcommand{\b}[1]{{\bf #1}} %
\newcommand{\Z}{\mathbb{Z}} %
\newtheorem{theorem}{Theorem} %
\newtheorem{lemma}{Lemma} %
\newtheorem{Qu}{Question}%
\newtheorem{conjecture}{Conjecture} %

\newcommand{\cp}{{\mathfrak p}} %

\begin{document}

\title{Artin's Conjecture on Zeros of $p$-Adic Forms}

\author{D.R. Heath-Brown}
\maketitle

\begin{abstract}
This is an exposition of work on Artin's Conjecture on the zeros of
$p$-adic forms. A variety of lines of attack are described, going back 
to 1945.  However there is particular emphasis on recent developments
concerning quartic forms on the one hand, and systems of quadratic
forms on the other.
\end{abstract}

Artin's Conjecture \cite[Preface]{Artin} is the following statement.
\begin{conjecture}
Let $F(x_1,\ldots,x_n)\in\Q_p[x_1,\ldots,x_n]$ be a form of degree $d$.
Then if $n>d^2$ the equation $F(x_1,\ldots,x_n)=0$ 
has a non-trivial solution in $\Q_p^n$.
\end{conjecture}
Here $\Q_p$ is the $p$-adic field corresponding to a rational prime 
$p$. Artin was led to his conjecture by considerations about
$C^i$-fields, and the above assertion can be re-phased to say that
$\Q_p$ is a $C^2$-field.
There are easy examples for every prime $p$ and every degree $d$ to show that
one cannot take $n=d^2$ here.  The conjecture can be generalized to
more general $\cp$-adic fields, and to
systems of forms of degrees $d_1,\ldots,d_r$, in which case the
condition on $n$ becomes $n>d_1^2+\ldots+d_r^2$.

One reason for the interest in Artin's Conjecture comes from the study
of Local-to-Global Principles.  One example is provided by the
following theorem of Birch \cite{B2}.
\begin{theorem}
Let $F(x_1,\ldots,x_n)\in\Q[x_1,\ldots,x_n]$ be a non-singular
form of degree $d\ge 2$ with $n>(d-1)2^d$.  Then
\begin{equation}
\label{birch}
\#\{\b{x}\in\Z^n:F(\b{x})=0,\;\max_{i=1,\ldots,n}|x_i|\le B\}
=c_FB^{n-d}+o(B^{n-d})
\end{equation}
as $B\rightarrow\infty$.  Moreover the constant $c_F$ is strictly
positive providing that the equation $F(x_1,\ldots,x_n)=0$ has zeros
in $\R^n$ and in each $p$-adic field.
\end{theorem}
Thus if Artin's Conjecture were true the $p$-adic condition would hold
automatically, since $(d-1)2^d\ge d^2$.

Unfortunately Artin's Conjecture is currently only know in the cases
$d=1$ and 2 (which are classical), and $d=3$ (due to Lewis \cite{L}).  
Indeed the conjecture is known to be false in general, the
first counterexample having been found by Terjanian 
\cite{T1}, for degree $d=4$. If one sets
\[G(x_1,x_2,x_3)=x_1^4+x_2^4+x_3^4-(x_1^2x_2^2+x_1^2x_3^2+x_2^2x_3^2)
-x_1x_2x_3(x_1+x_2+x_3)\]
and
\begin{eqnarray*}
F(x_1,\ldots,x_{18})&=&G(x_{1},x_{2},x_{3})+G(x_{4},x_{5},x_{6})+
G(x_{7},x_{8},x_{9})\\
&&\hspace{-1cm}\mbox{}+4G(x_{10},x_{11},x_{12})+4G(x_{13},x_{14},x_{15})
+4G(x_{16},x_{17},x_{18}),
\end{eqnarray*}
then $F(\b{x})$ is a form in 18 variables with 
only the trivial zero over $\Q_2$.  Subsequent
work has produced
counter-examples for many values of $d$, though $d$ is even in every
case known.

\begin{Qu}
Can one find any counter-examples to Artin's Conjecture with odd degree?
\end{Qu}

The most important general result in the positive direction is that of 
Ax and Kochen \cite{AK}.
\begin{theorem}
For every $d\in\N$ there is a $p_0(d)$ such
that Artin's Conjecture holds whenever $p\ge p_0(d)$.
\end{theorem}
The proof uses Mathematical Logic, and is based on the fact that the
analogue of Artin's Conjecture is known for the fields
$\mathbb{F}_p((t))$.  A value for $p_0(d)$ was found by Brown 
\cite{Brown}:-
\begin{equation}\label{bb}
2^{2^{2^{2^{2^{d^{11^{4d}}}}}}}!
\end{equation}
Here the ``!'' symbol is merely an exclamation mark, and not a
factorial sign! Another result by Ax and Kochen \cite{AK2}
shows that the theory of $p$-adic fields is decidable.
Thus for each fixed prime $p$ and each fixed degree $d$ there is, in
principle, a procedure for deciding whether the
statement\bigskip

``Every form $F(x_1,\ldots,x_{d^2+1})\in\Q_p[x_1,\ldots,x_{d^2+1}]$ has
a nontrivial zero 
\indent over $\Q_p$.''

\bigskip\noindent is true or false.
It follows that one can, in theory, test every prime up to Brown's
bound (\ref{bb}), and hence
decide whether or not Artin's Conjecture holds for a given degree $d$.

A second approach to Artin's Conjecture, developed by Lewis \cite{L}
for $d=3$, Birch and Lewis \cite{BL} for $d=5$, and Laxton and Lewis
\cite{LL} for $d=7$ and 11, applies a $p$-adic ``minimization'' process
to the form $F$ to produce a suitable model over $\Z_p$.  One then
examines the reduction $\overline{F}[\b{x}]\in\mathbb{F}_p[x_1,\ldots,x_n]$.
If this can be shown to have a non-singular zero, Hensel's Lemma will
allow us to lift it to a non-trivial zero of $F$ over $\Z_p$.  However
this ``minimization'' method has limited applicability. If 
$d$ can be written as a sum of composite numbers it is possible that 
$\overline{F}$ factors as $G_1^{e_1}\ldots G_k^{e_k}$ with ${\rm deg }\,
G_i\ge 2$ and $e_i\ge 2$ for every $i$.  In this case it is impossible
for $\overline{F}$ to have a non-singular zero.  The method is
therefore doomed to fail for such degrees.  In fact $d=1,2,3,5,7$ and
11 are the only integers which cannot be written as a sum of composite
numbers.  However for these values the method works moderately well, and
produces results of the type given by Ax and Kochen, but with much
smaller values for $p_0(d)$.  Thus Leep and Yeomans \cite {LY} 
showed that one may take $p_0(5)=47$, and Wooley \cite {W711}, that
$p_0(7)=887$ and $p_0(11)=8059$ are admissible.  These are susceptible
to further improvement, and indeed calculations by Heath-Brown have
shown that for $d=5$ Artin's Conjecture holds for $p\ge 17$.

\begin{Qu}
Does Artin's Conjecture hold for $d=5$, for every prime?
\end{Qu}
This is certainly decidable in principle, but whether it is realistic
to expect a computational answer with current technology is unclear.

The minimization approach can also be used for systems of forms.  It
shows (Demyanov \cite{Dem}) that $n>8$ suffices for a pair of 
quadratic forms, for every
$p$, and (Birch and Lewis \cite{BL2}, Schuur \cite{Sch}) that $n>12$ 
suffices for a system of 3
quadratic forms, providing that $p\ge 11$.  A very recent application
involving forms of differing degrees has been given by Zahid \cite{Z},
who shows that a quadratic and a cubic form over $\Q_p$ have a common
zero if $n>13=2^2+3^2$, providing that $p>293$.

Since Artin's Conjecture is false in general, it is natural to ask about the 
number $v_d(p)$, defined as the minimal integer such
that every form $F(x_1,\ldots,x_n)\in\Q_p[x_1,\ldots,x_n]$ of degree $d$
in $n>v_d(p)$ variables, has a non-trivial $p$-adic zero.  We also write 
$v_d=\max_p v_d(p)$.  Brauer \cite{Br} proved a result that implies that
$v_d$ is finite for every $d$.

\begin{theorem}
For every degree $d$ there is an integer $v_d$ such
that for each prime $p$, every form
$F(x_1,\ldots,x_n)\in\Q_p[x_1,\ldots,x_n]$ of degree $d$ with $n>v_d$ 
has a non-trivial $p$-adic zero.
\end{theorem}

Brauer's proof involves multiple nested inductions, and did not
lead to explicit bounds for $v_d$.  More recent versions of the
argument due to Leep and Schmidt \cite{Schm}, and particularly Wooley \cite{W}, 
are vastly more efficient, yielding
\begin{equation}\label{wb}
v_d\le d^{2^d}
\end{equation}
in general, but this is still disappointingly large.  Brauer's basic
idea is to show that for any $m\in\N$, the form $F$ will represent a 
diagonal form in $m$ variables as soon as $n$ is large enough compared
to $m$. It is not hard
to show (Davenport and Lewis \cite{DL})
that for every $p$ and every $d$ one can solve diagonal
equations
\[c_1x_1^d+\ldots+c_mx_m^d=0\]
over $\Q_p$ as soon as $m>d^2$.  Thus it suffices that $F$ should represent
a diagonal form in $m\ge d^2+1$ variables. We therefore seek
linearly independent vectors
$\b{e}_1,\ldots,\b{e}_m\in\Q_p^n$ such that
$F(\lambda_1\b{e}_1+\ldots+\lambda_m\b{e}_m)$ is a diagonal form in
$\lambda_1,\ldots,\lambda_m$.   If we choose the vectors $\b{e}_i$
inductively it is clear that $\b{e}_m$ must be a zero of a collection
of forms of degree strictly less than $d$.  Specifically there will be
$m-1$ forms of degree $d-1$;  $m(m-1)/2$ forms of degree $d-2$; and so
on.  The induction argument therefore involves the analogue of $v_d$
for systems of forms of differing degrees, and not just for a single
form of degree $d$.

There is an approach to these problems (Heath-Brown \cite{PLMS}) which is 
intermediate between the method of Lewis, Birch and Lewis, and Laxton
and Lewis and that of Brauer, Schmidt and Wooley. In this intermediate
approach one does not diagonalize $F$ fully, but removes enough of the
coefficients to ensure that there is a multiple of $F$ which has a
non-singular zero over $\mathbb{F}_p$, so that Hensel's Lemma can be used.
As an example we have the following lemma.
\begin{lemma}
Let $p\not=2,5$ or $13$ be prime and let
\[H(x,y,z)=Ax^4+Bxy^3+Cy^4+Dxy^3+Eyz^3+Fz^4\in\Q_p[x,y,z].\]
Suppose further that $A,C$ and $F$ are $p$-adic units. Then $H$
must represent zero non-trivially over $\Q_p$.
\end{lemma}
In order to produce such forms by the inductive construction above one
has to solve a system containing quadratic and linear equations, but
not cubics.

The power of this new method is well illustrated by the case $d=4$,
for which a direct application of (\ref{wb}) yields $v_4\le 4294967296$.
In contrast the new method (Heath-Brown
\cite[Theorem 2 and Note Added in Proof]{PLMS}) establishes the
following bounds.

\begin{theorem}\label{hb4th}
We have
\begin{enumerate}
\item[(i)] $v_4(p)\le 120$ for $p\ge 11$, 
\item[(ii)] $v_4(p)\le 128$ for $p=3$ and $p=7$, 
\item[(iii)] $v_4(5)\le 312$, and
\item[(iv)] $v_4(2)\le 4221$.  
\end{enumerate}
Thus $v_4\le 4221$
\end{theorem}

One sees that $p=2$ is the worst case by far.  It is fair to say that we have
absolutely no idea what the correct value for $v_4$ is, and it seems
natural in particular to ask the following question.

\begin{Qu}
Are there any counter-examples to Artin's Conjecture for quartic forms
with $p\not=2$?
\end{Qu}

It is convenient at this point to introduce the following notation.
For any field $K$, let
$\beta(r;K)$ be the least integer $m$ such that a system of $r$
quadratic forms over $K$ has a non-trivial common zero in $K$ as soon
as the number of variables exceeds $m$.  The case $d=2$ of Artin's
Conjecture, which is known to be true, yields $\beta(1;\Q_p)=4$, and
in general the conjecture would imply that $\beta(r;\Q_p)=4r$.  

The results on $v_4(p)$ from Heath-Brown \cite{PLMS} arise from the
estimates
\[v_4(p)\le\left\{\begin{array}{cc} 16+\beta(8;\Q_p), & p\not=2,5,\\
40+\beta(12;\Q_p), & p=5,\\  537+\beta(43;\Q_p), & p=2.\end{array}
\right.\]
together with suitable bounds for $\beta(r;\Q_p)$.  It is therefore
natural to turn our attention to the question of systems of
quadratic forms.  For general $r$ it has been shown by Leep
\cite{leep} that $\beta(r;\Q_p)\le 2r^2+2r$ for all $r$ and $p$.
There have been subsequent small improvements, but in all cases the bound is
asymptotic to $2r^2$ as $r\rightarrow\infty$.  Leep's argument is an
elementary induction on $r$, somewhat in the spirit of the Brauer
induction method.

A recent alternative attack (Heath-Brown \cite{Comp}) starts from
the work of Birch and Lewis \cite{BL2}, who used the minimization
approach to handle systems of three quadratic forms.  In general this
leads to a set of forms over $\mathbb{F}_p$ for which one wants to
find a non-singular common zero.  This is done via a counting
argument, so that one requires, amongst other information, an estimate
for the overall number of common zeros.  The following rather easy
lemma suffices.

\begin{lemma}\label{count2}
Suppose we have a system of quadratic forms 
\[Q^{(i)}(x_1,\ldots,x_n)\in \mathbb{F}_p[x_1,\ldots,x_n],\;\;\;(1\le i\le r)\]
with $N$ common zeros over $\mathbb{F}_p$.  
Write $N_R$ for the number of vectors
$\b{u}\in \mathbb{F}_p^r$ for which
\begin{equation}\label{comb}
\sum_{i=1}^r u_i Q^{(i)}(x_1,\ldots,x_n)
\end{equation}
has rank $R$, and assume that such a linear combination vanishes only
for $\b{u}=\b{0}$.  Then
\[|N-q^{n-I}|\le\sum_{1\le t\le n/2}q^{n-I-t}N_{2t}.\]
\end{lemma}

For vectors $\b{u}$ in the algebraic completion
$\overline{\mathbb{F}_p}$ the condition that (\ref{comb}) should have
rank at most $R$ defines a projective algebraic variety.  It is
possible to derive a good upper bound for the dimension of this set,
using the fact that the original $p$-adic system was minimized. This
bound on the dimension leads in turn to a bound for $N_R$.  This
enables one to show that the system of quadratic forms over
$\mathbb{F}_p$ has a non-singular zero when $p$ is large enough.  In
particular one can show that $\beta(r;\Q_p)=4r$ as soon as $p>(2r)^r$.

In contrast to the situation for the original formulation of Artin's
Conjecture, we know of no counter-examples for systems of quadratic
forms.  It is therefore possible that $\beta(r;\Q_p)=4r$ for every prime $p$.
\begin{Qu}
Is it true that $\beta(r;\Q_p)=4r$ for every prime $p$?
\end{Qu}
It is not even known what happens if we restrict the quadratic forms
to be diagonal.

The Ax-Kochen result already implies the existence of a bound $p_r$
such that $\beta(r;\Q_p)=4r$ for $p>p_r$. However the two methods have
a very important difference when we come to apply them to finite
extensions $\Q_{\cp}$ of $\Q_p$.  Suppose the residue field $F_{\cp}$ of such an
extension has cardinality $q=p^e$.  Then the Ax-Kochen theorem yields
the existence of a bound $p_{r,e}$ such that $\beta(r;\Q_{\cp})=4r$
for $p>p_{r,e}$.  Thus there is a condition on the characteristic of
$F_{\cp}$.  For example, the theorem leaves open the possibility that
$\beta(r;\Q_{\cp})>4r$ whenever $\Q_{\cp}$ is a finite extension of
$\Q_2$.  In contrast, the new method extends to give the following
result.
\begin{theorem}\label{hbthm}
We have $\beta(r;\Q_{\cp})=4r$ whenever $\#\mathbb{F}_{\cp}>(2r)^r$.
\end{theorem}
Here there is a condition on the cardinality of $\mathbb{F}_{\cp}$, rather
than its characteristic.

This makes a crucial difference when we consider the $u$-invariant of
function fields of the form $\Q_p(t_1,\ldots,t_k)$, as has been shown
by Leep \cite{Leepu}.  The $u$-invariant of a field $K$ is the
smallest integer $n$ such that any quadratic form over $K$ in more
than $n$ variables must have a non-trivial zero over $K$.  Thus
$u(\R)=\infty, u(\C)=1$ and $u(\Q_p)=4$. It is easy to see that
$u(K(x))\ge 2u(K)$ in general, and hence that
$u(\Q_p(t_1,\ldots,t_k))\ge 2^{2+k}$ for all $k\ge 0$.
Prior to the appearance of
the new results on $\beta(r;\Q_{\cp})$ just described, the only values
of $k$ for which it was known that $u(\Q_p(t_1,\ldots,t_k))$ is finite
were $k=0$ and $k=1$.  When
$k=1$, Parimala and Suresh \cite{PS} have recently shown that the 
$u$-invariant is 8, if $p\not=2$. Indeed Wooley, in work to appear,
has shown how to adapt the circle method to handle quite general
problems over $\Q_p(t)$, proving in particular that $u(\Q_p(t))=8$ for
every prime $p$.

In order to handle the $u$-invariant for function fields 
$\Q_p(\b{t})=\Q_p(t_1,\ldots,t_k)$ in $k$ variables, Leep considers 
a quadratic form 
$Q(X_1,\ldots,X_n)$ over $\Q_p(\b{t})$, in which the coefficients 
of $Q$ are polynomials in $t_1,\ldots,t_k$ of total degree at most
$d$, say.  One now considers a finite extension $\Q_{\cp}$ of $\Q_p$,
whose significance will become apparent later, and considers both $Q$
and the $X_i$ as polynomials in $t_1,\ldots,t_k$ over the new field
$\Q_{\cp}$. If we suppose that the 
$X_i$ are polynomials of total degree at most $D$ then the
overall number of coefficients in $X_1,\ldots,X_n$ is
\[N:=n(D+k)\ldots(D+1)/k!.\]
One may regard these coefficients as variables $c_1,\ldots,c_N\in \Q_{\cp}$, 
which one uses to force 
$Q(X_1,\ldots,X_n)$ to vanish identically.  Since $Q(X_1,\ldots,X_n)$
has total degree at most $2D+d$ as a function of $t_1,\ldots,t_k$
there are at most 
\[R:=(2D+d+k)\ldots(2D+d+1)/k!\]
coefficients which one must arrange to vanish.  Each of these is a
quadratic form in $c_1,\ldots,c_N$.  According to Theorem \ref{hbthm}
the corresponding system of quadratic forms has a non-trivial zero
$(c_1,\ldots,c_N)\in\Q_{\cp}$ providing that $N>4R$ and $q>(2R)^R$,
where $q$ is the cardinality of the residue field of
$\Q_{\cp}$. However it is clear that $N/R\rightarrow 2^{-k}n$ as $D\rightarrow
\infty$.  Hence if $n=1+2^{2+k}$ we can choose $D=D(k,d)$ so that
$N>4R$.  It follows that $Q(X_1,\ldots,X_n)=0$ has a non-trivial solution
$X_1,\ldots,X_n\in\Q_{\cp}(t_1,\ldots,t_k)$ providing that
$q>q_0(k,d)$.

One now calls on a result of Springer \cite{spr}, which states that if
$Q$ is a quadratic form over a field $F$ of characteristic different
from 2, which has a non-trivial zero
over some extension of $F$ of odd degree, then $Q$ has a non-trivial
zero over $F$ itself.  Thus to complete the proof it suffices to
choose $\Q_{\cp}$ to be an extension of $\Q$ of odd degree, and for
which $q>q_0(k,d)$.  One may then apply Springer's result with
$F=\Q_{\cp}(\b{t})$ to produce a non-trivial zero of $Q$ over the
original field $\Q_{p}(\b{t})$.  We therefore have the following
result, due to Leep \cite{Leepu}.
\begin{theorem}\label{lpthm}
We have $u(\Q_p(t_1,\ldots,t_k))=2^{2+k}$ for all
$k\in\N$ and all primes $p$.
\end{theorem}

The elegant feature of this argument is the way in which the size
constraint on $q$ disappears.  It is clear that the actual bound
$(2R)^R$ is irrelevant.  

One can utilise the case $k=1$ of Theorem \ref{lpthm} to obtain new
bounds for $\beta(r;\Q_p)$.  For example one has
$\beta(3;\Q_p)\le 16$ and $\beta(4;\Q_p)\le 24$ for every prime $p$.
These estimates are themselves used in the proof of Theorem
\ref{hb4th}.
It is curious that these results hold even for the case when the
residue field is small, even though Theorem \ref{hbthm}, from which
they derive, requires the residue field to be large.

As a corollary of Theorem \ref{lpthm} one can give an analogous
statement for pairs of quadratic forms.
\begin{theorem}\label{2forms}
Two quadratic forms over $\Q_p(t_1,\ldots,t_k)$, in at least
$1+2^{3+k}$ variables, have a non-trivial common zero.
\end{theorem}

This follows from a result of Brumer \cite {Brum}, which shows that if
$F$ is a field of characteristic different from 2, then a pair of
quadratic forms over $F$ will have a common zero as soon as the number
of variables exceeds $u(F(X))$.

As with Theorem \ref{lpthm}, there are examples showing that one
cannot reduce the number of variables.  Of course both results remain
true if we replace $\Q_p$ by a finite extension.

In conclusion we remark that it would be interesting to know what
happens for systems of cubic forms over $\Q_{\cp}$.  One might hope to
show that $r$ cubic forms in $n>9r$ variables have a common zero when
the cardinality $q$ of the residue field is large enough in terms of
$r$.  However this is currently known only for $r=1$, by the result of
Lewis \cite{L}.  If the general statement were established one could
deduce an analogue of Theorem \ref{lpthm} for cubic forms, with the
number of variables required to exceed $3^{2+k}$.  Here Springer's
theorem would be replaced by the observation that if $F$ is a field of
characteristic zero, then any cubic form with a zero over a quadratic
extension of $F$ also has a zero over $F$ itself.

\bigskip
\bigskip

Mathematical Institute, 

24--29, St  Giles', 

Oxford OX1 3LB, 

UK
\bigskip

{\tt rhb@maths.ox.ac.uk}

\end{document}